\documentclass[12pt]{article}

\usepackage{amsmath}
\usepackage{graphicx}
\usepackage{wrapfig}
\usepackage{a4wide}
\usepackage{caption}
\usepackage[T1]{fontenc}

\title{Topology and subsets - the story of a theorem}
\author{Micha{\l} Adamaszek}

\newcommand{\er}{\mathbf{R}}
\newcommand{\zet}{\mathbf{Z}}
\newcommand{\ce}{\mathbf{C}}

\newcommand{\sub}{\textrm{Sub}}
\date{}

\begin{document}

\maketitle

\footnotetext{Warwick Mathematics Institute and DIMAP, University of Warwick, Coventry, CV4 7AL, UK}
\footnotetext{This article is a translation of the original (in Polish), written for \textit{Matematyka, Spo{\l}ecze{\'n}stwo, Nauczanie, Vol.46}. It is an extended version of the talk by the author at the XLV Szko{\l}a Matematyki Pogl{\k{a}}dowej, Jachranka, Poland, 27-31 August 2010, \texttt{http://www.msn.ap.siedlce.pl}. The author thanks the editors for their kind permission to make the translation available to a wider audience and to Rupert Swarbrick for his comments on the first draft. }
\footnotetext{The author kindly acknowledges the support of the Centre for Discrete Mathematics and its Applications, EPSRC award EP/D063191/1, during the preparation of this text.}

This will be a story about topological spaces whose elements are subsets of other spaces. Let us first explain what that means.

Most of the spaces we would like to deal with arise from a number of basic building blocks, like the interval $I=[0,1]$ or the circle $S^1$. More complicated spaces are constructed using operations such as the Cartesian product $X\times Y$, the quotient $X/A$ or a more general quotient $X/\sim$, where we identify the points of $X$ in a way described by some relation $\sim$. For instance, gluing the endpoints of an interval yields a circle: $S^1=I\big/0\sim 1$. 

Sometimes we might need more complicated constructions. Suppose, for example, that we have $n$ identical gas particles, moving freely in a container of shape $X$, where $X$ is some space. We would like to describe the space of all possible locations of these particles. As a set, this space consists of all $n$-element subsets $C=\{x_1,\ldots,x_n\}$ of $X$. Now we say that the subsets obtained by a ``small'' movement of all particles in some fixed subset $C$ are ``close'' to $C$. The meaning of ``small'' depends on the topology of $X$ and, clearly, we only allow movements which do not make any two particles collide at the same position. This way we define small open neighbourhoods, therefore introducing the structure of a topological space in our family of subsets. Spaces of this kind are usually called configuration spaces. This particular construction can also be described formally as
$$C_n(X)=\{(x_1,\ldots,x_n)\in X^n: x_i\neq x_j\}\big/\Sigma_n,$$
where $\Sigma_n$ is the group of all permutations of the coordinates.

This example should have prepared the reader to appreciate the spaces which will be of our interest in this story:

\begin{center}
\begin{tabular}{ccp{6.5cm}}
$\sub_n(X)$ & $=$ & the space of \emph{at most} $n$-element, nonempty subsets of $X$.
\end{tabular}
\end{center}

Therefore a typical element of $\sub_n(X)$ is just the same as that of $C_n(X)$, namely a set of $n$ distinct points of $X$. Additionally, however, $\sub_n(X)$ contains all ``degenerate'' configurations, containing fewer than $n$ points. Small open neighbourhoods of a fixed subset $S$ are defined as previously, that is by ``small'' movements of all points in $S$, this time with the possibility of collision which results in a configuration of smaller cardinality. This also implies the reverse process where a point of a degenerate configuration can be replaced with a few distinct points very close to it.

\begin{wrapfigure}{l}{0.35\textwidth}
\includegraphics{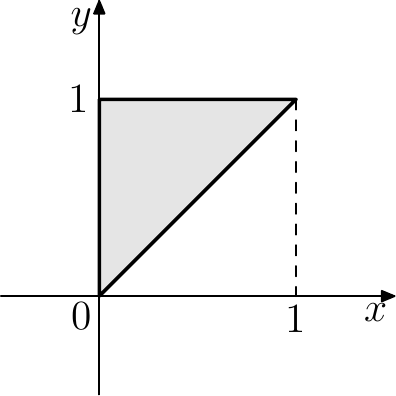}
\caption{$\sub_2(I)$}
\label{fig:1}
\end{wrapfigure}

We leave the formal definition of $\sub_n(X)$ as a topological space to the reader. Instead, let us illustrate this definition with an example. The space of at most two-element subsets of the interval $I$ can be represented as:
$$\sub_2(I)=\{(x,y): 0\leq x\leq y\leq 1\}$$
because every $1$- or $2$-element subset of the interval can be identified with a unique ordered pair. This situation is depicted in Fig.\ref{fig:1}. It follows that the space $\sub_2(I)$ is topologically equivalent (buzzword: homeomorphic) to a closed triangle. The points of the form $(x,x)$ correspond to singletons $\{x\}$. We see that a small neighbourhood of, say, the singleton $\{\frac{1}{2}\}$, contains subsets of the form $\{x,y\}$, where $x$ and $y$ are very close to $\frac{1}{2}$ and not necessarily equal.

Before proceeding further let us note two simple properties, which hold for an arbitrary space $X$:
$$\sub_1(X)=X$$
$$X=\sub_1(X)\subset \sub_2(X)\subset \sub_3(X)\subset\cdots$$
It is now time to explain the goal of this article, which is to understand the space
$$\sub_3(S^1)$$
consisting of at most $3$-element subsets of the circle. In particular, is it topologically equivalent to some well-known space?

\begin{wrapfigure}{l}{0.25\textwidth}
\includegraphics{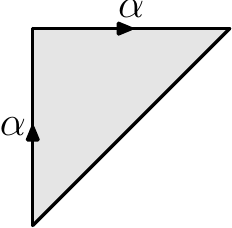}
\caption{$\sub_2(S^1)$}
\label{fig:2}
\end{wrapfigure}

As a warm-up we shall look at the space $\sub_2(S^1)$. Because the circle is the interval with identified endpoints ($S^1=I/0\sim 1$), also the space $\sub_2(S^1)$ arises from $\sub_2(I)$ via the identification of all occurrences of $0$ and $1$. It means that we have
$$\sub_2(S^1)=\sub_2(I)\big/(0,x)\sim (x,1).$$

\begin{figure}[h]
\begin{center}
\begin{tabular}{cccc}
\includegraphics{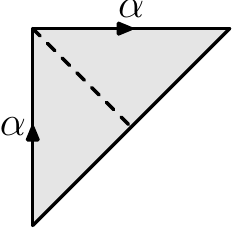} & \includegraphics{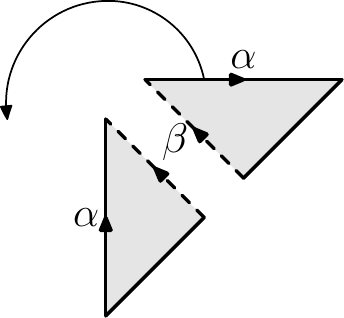} & \includegraphics{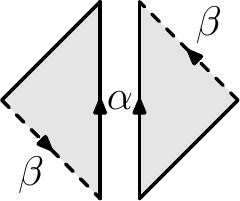} & \includegraphics{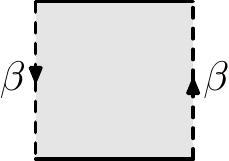} 
\end{tabular}
\end{center}
\caption{$\sub_2(S^1)$ is the M\"obius band}
\label{fig:3}
\end{figure}

Let us see what happens to the triangle of Fig.\ref{fig:1} under this identification. The left edge contains the points $(0,x)$ and the top edge contains the points $(x,1)$, so we must glue those two edges along the arrows in Fig.\ref{fig:2}. One can construct a model out of paper, looping the left edge into a circle and winding the top edge onto it (note that after the identification all three vertices $(0,0)$, $(0,1)$ and $(1,1)$ represent the same singleton $\{0\}$). Since it may not be immediately clear what we get, let us adopt a different approach, outlined in Fig.\ref{fig:3}. First we cut the triangle along the dashed line and move the two pieces apart, remembering about the identification along the dashed arrows now called $\beta$. Next we glue the triangles along the edge $\alpha$ and we easily recognize the M\"obius band --- a space obtained from a strip of paper by the identification of one pair of opposite edges with opposite orientations. We have proved the following theorem.

\textbf{Theorem.}\textit{
The space $\sub_2(S^1)$ of at most two-element subsets of the circle is topologically equivalent to the M\"obius band. It contains the space $\sub_1(S^1)$ of one-element subsets as the boundary circle.}

To obtain the second statement one just needs to track the points $(x,x)$ during the cut-and-glue process.

To entertain all more demanding readers we will now sketch a different argument. Every two points on a circle determine a line (with a single point we associate the line tangent to the circle at that point). The direction of that line gives an element of the one-dimensional projective space $\er P^1$, topologically equivalent to $S^1$. This way we defined a continuous map $\sub_2(S^1)\to S^1$. The fibre of this map (the pre-image of any point) consists of all lines of a fixed direction which intersect the circle, so it is equivalent to the closed interval. It therefore follows that $\sub_2(S^1)$ is topologically equivalent to a bundle of intervals over a circle. There are just two such bundles: the cylinder $S^1\times I$ and the M\"obius band. The reader who got this far can easily eliminate the first possibility.

After this warm-up let us look at spaces of three-element subsets, starting with the interval and $\sub_3(I)$. Every 3-, 2- or 1-element subset of the interval can be written as an ordered triple, but the 2-element subsets
$\{x,z\}$ have two such representations $(x,x,z)$ and $(x,z,z)$, which must be identified. It follows that
$$\sub_3(I)=\{(x,y,z): 0\leq x\leq y \leq z \leq 1\}\big/(x,x,z)\sim (x,z,z).$$
\newpage 

\begin{wrapfigure}{l}{0.4\textwidth}
\includegraphics{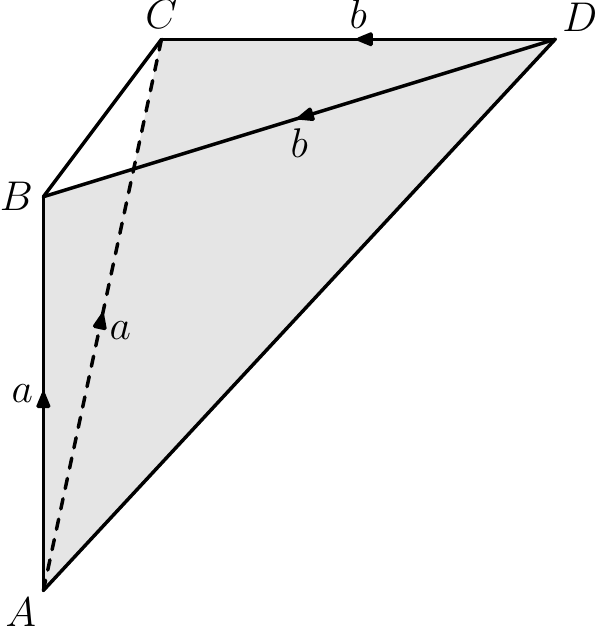}
\caption{$\sub_3(I)$}
\label{fig:7}
\end{wrapfigure}

The set $\{(x,y,z): 0\leq x\leq y \leq z \leq 1\}$ is a tetrahedron with vertices $A=(0,0,0)$, $B=(0,0,1)$, $C=(0,1,1)$, $D=(1,1,1)$. The space $\sub_3(I)$ arises from this tetrahedron by gluing the face $ABD$ (points of the form $(x,x,z)$) to $ACD$ (points of the form $(x,z,z)$) in such a way that $B$ is identified with $C$, $AB$ with $AC$ and $BD$ with $CD$, as in Fig.\ref{fig:7}. The final product is easy to visualize if we imagine the tetrahedron is stretchable and we rotate $C$ around the edge $AD$ until the shaded faces coincide. We get a solid which looks like two filled cones with apexes $A$ and $D$, glued along the bases. Topologically this is just a closed ball. The subspace of one-element subsets of the form $(x,x,x)$ is contained inside as the diameter $AD$. By the way, more advanced readers can deduce (how?) that our target space $\sub_3(S^1)$ is a $3$-dimensional manifold (meaning that each point has a neighbourhood which is an open $3$-dimensional ball).

The space $\sub_3(S^1)$ can now be obtained, as before, by an identification of $0$ and $1$. Therefore:
$$\sub_3(S^1)=\sub_3(I)\big/(0,y,z)\sim (y,z,1).$$

\begin{wrapfigure}{r}{0.4\textwidth}
\includegraphics{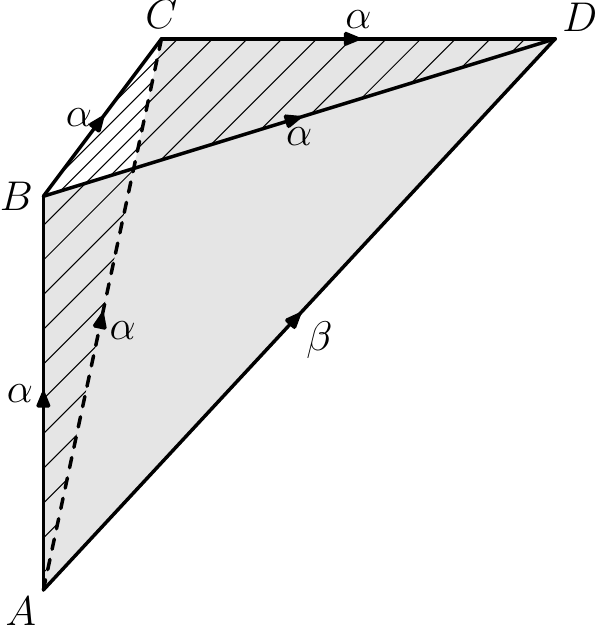}
\caption{$\sub_3(S^1)$}
\label{fig:8}
\end{wrapfigure}

As a result, the way to obtain $\sub_3(S^1)$ from the tetrahedron $ABCD$ is to glue the faces $ABD$ and $ACD$ (as in $\sub_3(I)$) and, additionally, glue the faces $ABC$ (the points $(0,y,z)$) and $BCD$ (the points $(y,z,1)$) so that the edges $AB$, $BC$ and $CA$ are identified, respectively, with $BC$, $CD$ and $DB$. These identifications are indicated in Fig.\ref{fig:8}. The latter one seriously complicates the issue. Our previous techniques --- cutting, gluing or building models --- may not be sufficient to handle this case.

In fact the first ones to describe the topology of $\sub_3(S^1)$ were Karol Borsuk (1905-1982) and Raoul Bott (1923-2005). Karol Borsuk is known in contemporary mathematics, among other things, for the definition of a cofibration, for the Borsuk Conjecture (Can every convex set in $R^d$ be divided into $d+1$ sets of smaller diameter?), disproved in 1993 and for the Borsuk-Ulam theorem, whose popular formulation states that at each moment there exist two antipodal points on the surface of the Earth with the same temperature and the same atmospheric pressure. The second one, Raoul Bott, a mathematician of Hungarian origin, is associated mainly with the famous, initially very surprising periodicity theorem for the unitary group, which soon after being proved became one of the fundamental pieces of the emerging $K$-theory.
 
The story of the answer, given by Borsuk and Bott, is just as interesting as the answer itself. In 1949 Borsuk published his paper \cite{Borsuk}, where he proved that $\sub_3(S^1)$ is topologically equivalent to the product $S^2\times S^1$ of the standard sphere and a circle. Three years later Fundamenta Mathemathicae published a paper by Bott \cite{Bott}. It has the form of a letter to Borsuk, starting with the words (original notation):

\begin{center}
\parbox{0.8\textwidth}{\small
In your paper \textit{On the third symmetric potency of the circumference} (...) you assert that the third symmetric potency of $S_1^{(3)}$ of the circle $S_1$ is homeomorphic to the Cartesian product of $S_1$ and the two sphere $S_2$. (...). But in fact the identification you have made is incorrect and in consequence your final conclusion (...) is false. A quite simple and short argument shows that $S_1^{(3)}$ has a vanishing fundamental group whence (...) $S_1^{(3)}$ is a simply connected lensespace, \textit{i.e.} the three sphere $S_3$.
}
\end{center}

It means that Borsuk made a mistake and the correct answer should have been:

\textbf{Theorem.}\textit{The space $\sub_3(S^1)$ is topologically equivalent to the three-dimensional sphere $S^3$.}

If the reader shivers at the prospect of the $3$-sphere, the easiest way around is to think about it as the $3$-dimensional space $\er^3$ with an additional ``point at infinity''. This is an analogy with similar descriptions of the lower-dimensional spheres.

At this point both the reader and Karol Borsuk deserve a few words of an explanation. Borsuk's paper, although quite complicated, is mostly correct. The mistake crept in only at the very end of the argument. This fact is also emphasized by Bott, who corrects only the conclusion of Borsuk's paper. It goes as follows. After some rather long manipulation Borsuk proves that $\sub_3(S^1)$ can be obtained from two solid tori via some identification of the boundaries. The confusing bit is that there are two possible identifications: the meridians of one torus can be glued with the meridians or with the parallels of the other one. In the first case we obtain $S^2\times S^1$ (quite easy to see), while the other case yields $S^3$ (a bit harder).

The readers familiar with some basic magic spells of algebraic topology can themselves decide who is right. Our description of $\sub_3(S^1)$ as the quotient space of the tetrahedron $ABCD$, although difficult to imagine, is a perfectly valid cell decomposition (strictly speaking: a $\Delta$-complex decomposition). It is therefore suitable for the calculation of topological invariants such as the fundamental group and homology groups. Let us compute, following Bott, the fundamental group $\pi_1(\sub_3(S^1))$ of $\sub_3(S^1)$. We begin with a quick reminder (or introduction) of the definition. The fundamental group of a space $X$ is generated by all loops (continuous functions $f:S^1\to X$) with a fixed basepoint $x_0$ (which means $f(1)=x_0$). We identify two loops $f$ and $g$ if one of them can be continuously deformed into the other, by which we mean the existence of an entire family of loops varying continuously from $f$ to $g$. The technical term for such two loops is ``homotopic''. The elements of the fundamental group $\pi_1(X)$ are the equivalence classes of loops under the homotopy relation.

All this sounds like it is very technical, but with the cell decomposition of a space at hand it is easy to describe its fundamental group. First of all we can restrict our attention to the loops in the one-dimensional skeleton. In our space there are two such loops, $\alpha$ and $\beta$ (recall that all the points $A$, $B$, $C$ and $D$ were identified). Each two-dimensional face bounds some loop, thus trivializing it --- such a loop can be contracted to a point (a trivial loop) via a deformation which pulls it across the face. We see from Fig.\ref{fig:8} that the face $ABD$ (as well as $ACD$) has as its boundary the loop $\alpha\alpha\beta^{-1}=\alpha^2\beta^{-1}$, while the face $ABC$ (and $BCD$) bounds the loop $\alpha\alpha\alpha^{-1}=\alpha$. It follows that the group $\pi_1(\sub_3(S^1))$ is generated by two elements $\alpha$, $\beta$ with the relations $\alpha^2\beta^{-1}=1$ and $\alpha=1$. We use the notation
$$\pi_1(\sub_3(S^1))=\langle\alpha, \beta\ |\ \alpha=1, \alpha^2\beta^{-1}=1\rangle.$$
All this means that $\alpha=1$ and $\beta=\alpha^2=1$, so the group in question is trivial. This itself suffices to reject the answer $S^2\times S^1$, because the fundamental group of the latter space is $\zet$. Moreover, the readers who are up-to-date with the latest developments in mathematics and who checked a few paragraphs ago that $\sub_3(S^1)$ is a manifold can now conclude that it must be the sphere $S^3$. It is a corollary of the Poincar\'e Conjecture in dimension $3$, one of the Millennium Problems of the Clay Institute, proved recently by Grisha Perelman. It states precisely that $S^3$ is the only compact three-dimensional manifold with trivial fundamental group. We have verified that Bott's correction was fully justified.

Having sorted this out, let us move on to another intriguing question. We already know that $\sub_3(S^1)$ is the sphere $S^3$, or, in other words, the one-point compactification of the three-dimensional space $\er^3$. It contains the subspace $\sub_1(S^1)$ of one-element subsets of $S^1$. However, we know that $\sub_1(S^1)=S^1$ is a circle and a circle in $\er^3$ is usually called a knot. Which kind of a knot is it?

\begin{wrapfigure}{r}{0.4\textwidth}
\includegraphics[scale=0.7]{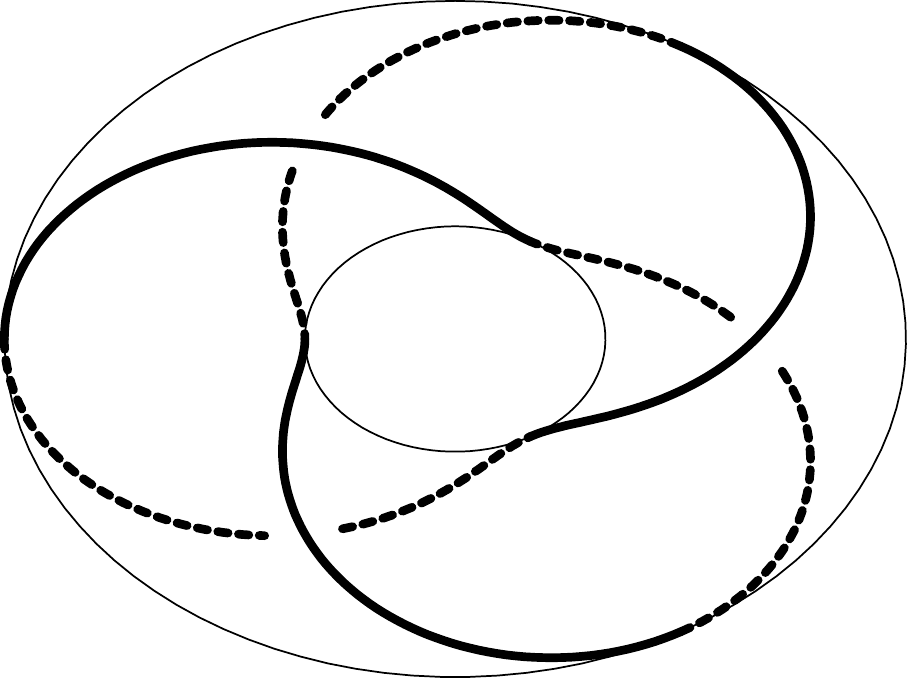}
\caption{Trefoil knot}
\label{fig:9}
\end{wrapfigure}

This question is answered for instance in a new paper by Jacob Mostovoy \cite{Mostovoy}. To encourage the reader to read that short note we will only mention that Mostovoy provides an explicit, analytic formula for the continuous bijection between $\sub_3(S^1)$ and $S^3$. With such a formula one can derive the equation satisfied by the subspace $\sub_1(S^1)$. We will just quote the final answer. If we identify the sphere $S^3$ with the set $\{(u,w)\in\ce\times\ce: |u|^2+|w|^2=1\}$, then the subspace $\sub_1(S^1)$ is given by the equation $u^3=w^2$. After some rescaling it can also be written in a parametric form
$$S^1\ni z\to (z^2,z^3)\in S^1\times S^1\subset S^3$$
which shows that our knot lies on the torus $S^1\times S^1$ and circles around it twice in the parallel direction and three times in the meridian direction, as in Fig.\ref{fig:9}. Such a knot is called a $(2,3)$-torus knot (the $(p,q)$-torus knots are defined in a similar fashion). This specific knot is the simplest of nontrivial knots, called the \emph{trefoil knot}.

We will not reproduce Mostovoy's argument here, but we encourage the readers to analyze it on their own. Instead we will present another proof of the trefoil knot theorem, using what we have already learned about the fundamental group. Our strategy is to calculate the \emph{knot group}. It is an invariant defined for any knot $K$ as
$$\pi_1(S^3\setminus K)$$
which is the fundamental group of what remains from $\er^3$ after removing the knot. In our problem the sphere $S^3$ is represented by the tetrahedron $ABCD$ with appropriate identifications on the boundary (Fig.\ref{fig:8}), and the knot $K$ is contained inside as the set of points of the form $(x,x,x)$, which is the segment $\beta=AD$. In this model the removal of the knot $K$ is equivalent to the removal of $\beta$ and the vertices $B$ and $C$ (which represent the same point as $A$ and $D$). 

\begin{wrapfigure}{l}{0.45\textwidth}
\includegraphics{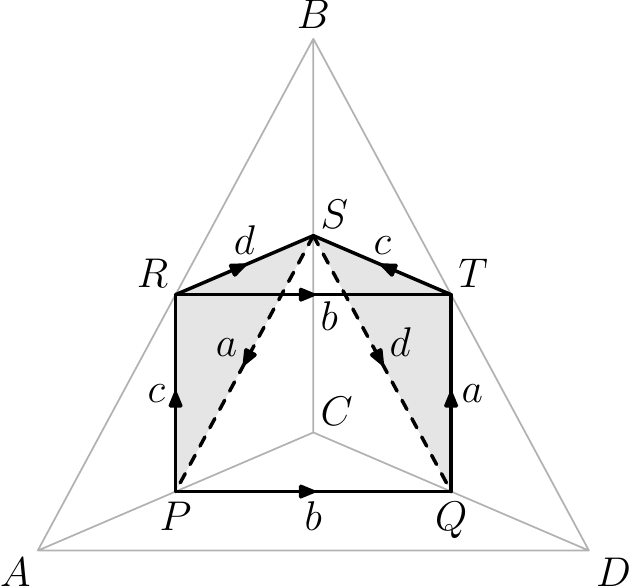}
\caption{}
\label{fig:10}
\end{wrapfigure}

To obtain a nicely triangulated space it is more convenient to remove the knot $K$ together with a small open neighbourhood. In our model this step is effected by removing a narrow prism along the segment $\beta$ together with small tetrahedral caps around $B$ and $C$. Finally, we can employ a continuous deformation, which does not influence the fundamental group, to enlarge the removed neighbourhoods until each of them reaches the middle of the edge $\alpha$. What remains is the space depicted in Fig.\ref{fig:10} --- a pyramid with base  $PRTQ$ and apex $S$. The identifications in $ABCD$ translate to identifications of faces and edges of that pyramid indicated in the figure (please check!). The one-dimensional cells $a,b,c,d$ are the loops which generate the fundamental group $\pi_1(S^3\setminus K)$. As before, we are going to write out the relations which hold in that group:
\begin{itemize}
\item top face $RTS$: $bcd^{-1}=1$,
\item front face $QTRP$: $ab^{-1}c^{-1}b=1$,
\item left face $SPR$ (also right face $QTS$): $acd=1$,
\item back face $SPQ$: $abd^{-1}=1$.
\end{itemize}
The first relation implies $c=b^{-1}d$ and the second one gives $a=b^{-1}cb=b^{-1}b^{-1}db=b^{-2}db$. The third relation takes the form
$$1=acd=b^{-2}dbb^{-1}dd=b^{-2}d^3,$$
so we have $b^2=d^3$. The last relation turns out to be redundant:
$$1=abd^{-1}=b^{-2}dbbd^{-1}=b^{-2}db^2d^{-1}=b^{-2}dd^3d^{-1}=b^{-2}d^3=1.$$
Eventually, the fundamental group we are computing has two generators $b$, $d$ and one relation $b^2=d^3$:
$$\pi_1(\sub_3(S^1)\setminus \sub_1(S^1)) = \langle b,d\ |\ b^2=d^3\rangle.$$
This group is quite well known. It is the so-called braid group on $3$ strands. The initial chapters of any textbook on knot theory identify it as the group of the trefoil knot (more generally, the group of any $(p,q)$-torus knot supports a single relation $x^p=y^q$). Does it imply that the knot under consideration must be the trefoil knot? In general, this question is rather delicate, because two different knots can have the same groups (explicit examples of such pairs are known). This time, however, we are lucky, because the trefoil knot, as well as any $(p,q)$-torus knot, is uniquely determined by its group, which in this full generality is a result of Burde and Zieschang. Our trefoil knot theorem is therefore proved.

This is the end of the story. Along the way we encountered the M\"obius band, we learned to calculate the fundamental group of a space, we constructed an unusual model of the sphere $S^3$, we applied the Poincar\'e Conjecture, we drew knots on a torus, we learned how to compute a group-theoretic knot invariant and all this was an outcome of a seemingly inconspicuous problem. Let us summarize our achievements in a theorem.

\textbf{Theorem (Borsuk-Bott-Reader-Mostovoy).}  \textit{The space $\sub_3(S^1)$ of at most three-element subsets of the circle is topologically equivalent to the $3$-sphere $S^3$. The subspace $\sub_1(S^1)$ is contained in it as the trefoil knot and the space $\sub_2(S^1)$ is a M\"obius band, whose boundary is that knot.}

Now, what does a M\"obius band with the trefoil knot on the boundary look like?

\end{document}